\documentclass[12pt]{amsart}

\usepackage[top=4cm,bottom=4cm,inner=4cm,outer=4cm]{geometry}
\usepackage{comment}
\newcommand{\disp}{\displaystyle}
\newcommand{\nc}{\newcommand}
\nc{\G}{{\Gamma}} \nc{\BC}{{\mathbb C}} \nc{\BQ}{{\mathbb Q}}
\nc{\BR}{{\mathbb R}} \nc{\BZ}{{\mathbb Z}} \nc{\BP}{{\mathbb P}} \nc{\PC}{{\BP_1(\BC)}}
\nc{\BN}{{\mathbb N}} \nc{\BM}{{\mathbb M}}
\nc{\fH}{{\mathbb H}}

\nc{\mat}{{\binom{a\,\ b}{c\,\ d}}}
\nc{\U}{{\mathcal U}}
\nc{\PS}{{\mbox{PSL}_2(\BZ)}} \nc{\SL}{{\mbox{SL}_2(\BZ)}}
\nc{\SR}{{\mbox{SL}_2(\BR)}} \nc{\PR}{{\mbox{PSL}_2(\BR)}}
\nc{\SLC}{{\mbox{SL}_2(\BC)}}

\nc{\GL}{{\mbox{GL}}} \nc{\PQ}{{\mbox{PGL}_2^+(\BQ)}}
\nc{\GR}{{\mbox{GL}_2^+(\BR)}} \nc{\PG}{{\mbox{PGL}_2(\BC)}}
\nc{\GC}{{\mbox{GL}_2(\BC)}}
\nc{\f}{{\mathcal{F}(\fH)}}
\nc{\Cc}{\widehat{\BC}}
\nc{\e}{{E_{\rho}(\G)}}
\nc{\g}{{\gamma}}
\nc{\vm}{{V_{\rho}(\G)}}
\nc{\oo}{{\mathcal O}}
\nc{\M}{{\mbox{M}}}
\nc{\om}{{\omega}}
\nc{\Om}{{\Omega}}
\nc{\TX}{{\widetilde{X}}}
\nc{\ol}{\overline}
\nc{\cl}{{\mathcal L}}
\nc{\ce}{{\mathcal E}}
\nc{\la}{{\lambda}}
\nc{\La}{{\Lambda}}

\nc{\cz}{{\mathcal Z}}

\newtheorem{numbered}{}[section]
\newtheorem{thm}[numbered]{Theorem}

\newtheorem{remark}[numbered]{Remark}
\newtheorem{prop}[numbered]{Proposition}

\newtheorem{example}[numbered]{Example}
\numberwithin{equation}{section}

\newcommand{\propref}[1]{Proposition~\ref{#1}}

\begin{document}
	
	\title[]{Modular differential equations and algebraic systems}
	\author[]{Hicham Saber} \author[]{Abdellah Sebbar}
	\address{Department of Mathematics, Faculty of Science, University of Ha'il,   Ha'il, Kingdom of Saudi Arabia}
	\address{Department of Mathematics and Statistics, University of Ottawa,
		Ottawa Ontario K1N 6N5 Canada}
	\email{hicham.saber7@gmail.com}
	\email{asebbar@uottawa.ca}
\maketitle

	\begin{abstract}
		In this paper, we show how solutions to  explicit algebraic systems lead to solutions to infinite families of modular differential equations.
	\end{abstract}

	\section{Introduction}

The theory of modular differential equations, which are linear differential equations with coefficients in the ring of modular forms, have been considered by early automorphic forms experts such as  Klein \cite{klein}, Hurwitz \cite{hurwitz} and Van der Pol \cite{vdp}. There has been a lot of interest in these differential equations in recent decades starting with the pioneering work by Kaneko and Zagier \cite{ka-za}. The subject developed into a fertile research area with applications in many areas of mathematics  and mathematical physics. A great deal of literature has been written on the subject, including the works \cite{yang,fr-ma,grabner,ka-ko,ka-et-al, milas,mukhi,nakaya}. We shall be concerned with modular differential equations in connection with the Schwarz differential equation and the theory of equivariant functions as we now explain.

 Let $\G$ be a discrete subgroup of $\SR$ acting on the upper half-plane $\fH$, and denote by $\ol{\G}$ its projection in $\PR$. We consider the following differential equation with an automorphic potential
\[
y''\,+\,Q(\tau)\,y\,=\,0,\ \ \tau\in\fH,
\]
where $Q(\tau)$ is a weight 4 automorphic form for $\G$. If $f_1$ and $f_2$ are linearly independent solutions, then $h=f_2/f_1$ satisfy the Schwarz differential equation
\[
\{h,\tau\}\,=\,2Q(\tau),
\]
where $\{h,\tau\}$ is the Schwarz derivative defined by
\[
\{h,\tau\}\,=\, \left(\frac{h''(\tau)}{h'(\tau)}\right)'\,-\,\frac{1}{2}\left(\frac{h''(\tau)}{h'(\tau)}\right)^2.
\]
The Schwarz derivative has many projective, geometric and analytic properties that can be found in \cite{nehari,s-g}.
On the other hand, for a meromorphic function $h$ on $\fH$, one can show that $\{h,\tau\}$ is a weight 4 automorphic form for $\G$ if and only if there exists a 2-dimensional representation $\varrho$ of $\ol{\G}$ such that
\[
h(\gamma\cdot\tau)\,=\,\varrho(\gamma)\cdot h(\tau)\,,\ \tau\in\fH\,,\ \gamma\in\G.
\]
We call  such a function $h$ a $\varrho-$equivariant function for $\G$. This class of functions has been studied extensively in \cite{structure,rational,vvmf,kyushu} with interesting applications in \cite{critical, forum, ramanujan,jmaa,baus}. 
The automorphic functions (of weight zero) are $\varrho-$equivariant with $\varrho=1$; the constant representation. 
If $\varrho=Id$, the defining representation, then $h$ is simply called an equivariant function (it commutes with action of $\G$). As an example, if
$f$ is a weight $k$ automorphic form for $\G$, then
\[
h_f(\tau)\,=\, \tau\,+\,k\,\frac{f(\tau)}{f'(\tau)}
\]
is an equivariant function for $\G$. This also includes  the case $f$ being a non-constant automorphic function which leads to the trivial equivariant function $h(\tau)=\tau$.  

In this paper, we focus on the case of the modular group $\G=\SL$. A holomorphic weight 4 modular form $Q(\tau)$ is thus a scalar multiple of the weight 4 Eisenstein series $E_4(\tau)$. Therefore, we consider
the modular differential equation
\begin{equation}\label{equ1}
y''\,+\,s\,E_4(\tau)\,y\,=\,0,
\end{equation}
and the corresponding Schwarz differential equation
\begin{equation}\label{equ2}
\{h,\tau\}\,=\,2\,s\,E_4(\tau)\,.
\end{equation}
It should be noted that the modular differential equations studied in \cite{ka-ko} and \cite{ka-za} can be reduced to the equation \eqref{equ1}, \cite{baus}.
According to \cite{forum}, any solution $h$ to \eqref{equ2} is necessarily locally univalent and leads to solutions $y_1=1/\sqrt{h'}$ and $y_2=h/\sqrt{h'}$ to \eqref{equ1}.
Moreover, for a solution to \eqref{equ2} to be meromorphic or to have a logarithmic singularity at $\infty$, the parameter $s$ must satisfy $s=\pi^2r^2$ where $r$ is a rational number.

In \cite{forum}, we investigated solutions to \eqref{equ2} that are $\varrho-$equivariant with $\mbox{Ker}\,\varrho$ having a finite index in $\SL$, in other words, that are modular functions for a  finite index subgroup of $\SL$. It turns out that necessarily $\rho$ is an irreducible representation of $\SL$ and that $s=\pi^2r^2$ with a rational number $r=n/m$, $2\leq m\leq 5$ and
$\gcd(m,n)=1$. Furthermore, the solution $h$ is a modular function for the principal congruence subgroup $\G(m)$. The integers $m$ and $n$ have the following interpretation:
 We have the two coverings of compact Riemann surfaces 
 \[
\pi:X(\ker\varrho)\longrightarrow X(\SL)\cong \PC
\]
 induced by the natural inclusion $\ker\varrho\subseteq\SL$,
and
\[
h:X(\ker\varrho)\longrightarrow X(\SL)\cong \PC
\] induced by the solution $h$. 
 Here $X(\G)$ is the modular curve attached to the subgroup $\G$. Then $m$ and $n$ are the respective ramification indices above $\infty$ for the two coverings.

 When $r$ is an integer $(m=1)$, the situation is completely different. There are always solutions to \eqref{equ2} that are simply equivariant, that is,  when $\varrho=Id$,  while the solutions to \eqref{equ1} are constructed from quasi-modular forms  \cite{jmaa}.
 
 In \cite{ramanujan}, we investigated the case when solutions to \eqref{equ2} correspond to reducible representation $\varrho$ of $\SL$. It turns out that necessarily $r=n/6$ with $\gcd(n,6)=1$. The denominator 6 occurs because it is the level of the commutator group of $\SL$ over which the  characters of $\SL$ are trivial.
 In addition, the solutions to \eqref{equ2} are integrals of weight 2 meromorphic modular forms with a character. For the case $n=1$, the weight 2 form in question is $\eta^4$. We then constructed solutions for every $n=1+12k$, $k\in\BN$, by integrating the modular form
 \[
 f_n=\frac{\eta^4}{\prod_{i=1}^{k}(J-a_i)^2},
 \]
 where the numbers $a_i$ are solution to the algebraic system 
\[
\frac{4}{x_i}\,+\,\frac{3}{x_i-1}\,+\,\sum_{j\neq i}\,\frac{12}{x_j-x_i}\,=\, 0\,,\ \ 1\leq i\leq k,
\]
which turns out to admit a solution in $(0,1)^k$.
The idea is to adjoin double poles to $\eta^4$ in $\fH$ with zero residues. In this case, the double poles are not elliptic points.

In this paper, we show that this method extends nicely to the remaining cases of residues modulo 12, namely, for $n$ coprime with 6 such that $n\equiv 5\,,7\,\mbox{or}\,11\mod 12$.
More precisely, starting from a fundamental solution $f_{\alpha}$ to 
\eqref{equ2} with $s=\pi^2(\alpha/6)^2$ for each  $\alpha=5$, 7 or 11,  one can construct
a solution for each $m$ in the residue class of $\alpha$ modulo 12 by adjoining double poles to $f_{\alpha}'$ and integrating. In this cases, the double poles are allowed to include one of the elliptic points $i$ or $\rho$ or both. However, it is shown that the whole construction can be carried out by solving the algebraic system 
\[
\frac{a}{x_i}\,+\,\frac{b}{x_i-1}\,+\,\sum_{j\neq i}\,\frac{c}{x_j-x_i}\,=\, 0\,,\ \ 1\leq i\leq k,
\]
where $a$, $b$ and $c$ vary with $\alpha$.

Furthermore, we revisit the cases where the level $m\in\{2,3,4,5\}$ studied in \cite{forum} and we show that our method can be applied to construct the solutions to \eqref{equ1} and \eqref{equ2}. Indeed, starting from a solution $t_m$ corresponding to $r=1/m$, which turns out to be a Hauptmodul for $\G(m)$, we can construct a solution corresponding to $r=n/m$ with $n$ coprime with $m$ by adjoining to $t_m'$ double zeros and double poles that arise from solutions to the above algebraic system with an appropriate choice of the parameters $a$, $b$ and $c$.

It is yet to be fully understood why solutions to a simple algebraic system would lead to solutions to infinite families of modular differential equations.

\section{Special values of higher derivatives of modular forms}
In this section we recall some classical elliptic modular forms. We also establish some interesting identities involving special values of their higher derivatives at elliptic fixed points. 

The Eisenstein series $E_2$, $E_4$ and $E_6$ are defined by their $q-$expansions:
\begin{align*}
E_2(\tau)&=1-24\,\sum_{n\geq 1}\,\sigma_1(n)\,q^n\,,\\
E_4(\tau)&=1+240\,\sum_{n\geq 1}\,\sigma_3(n)\,q^n\,,\\
E_6(\tau)&=1-504\,\sum_{n\geq 1}\,\sigma_5(n)\,q^n\,.
\end{align*}
Here $\tau$ is a variable in the upper half-plane $\fH=\{\tau\in\BC|\mbox{Im}(\tau)>0\}$ and  $q=\exp(2\pi i\tau)$ is the uniformizer at $\infty$. The arithmetical function $\sigma_k$ is defined on positive integers by
\[
\sigma_k(n)=\sum_{0<d\,|\,n}\,d^k.
\]
The function $E_2$ is a quasi-modular form of weight 2 and $E_4$ and $E_6$ are modular forms of respective weights 4 and 6 for the full modular group $\SL$.

We also define the Dedekind eta-function by
\[
\eta(\tau)\,=\,q^{\frac{1}{24}}\,\prod_{n\geq1}(1-q^n)\,,
\]
and the weight 12 cusp form $\Delta$ (the modular discriminant)
\[
\Delta(\tau)\,=\,\eta(\tau)^{24}\,=\,\frac{1}{1728}(E_4(\tau)^3-E_6(\tau)^2).
\]
We also have the elliptic modular function $J$  ($J-$invariant)
\[
J(\tau)\,=\,\frac{1}{1728}\frac{E_4(\tau)^3}{\Delta},
\]
and the Klein elliptic modular function $\lambda$ for $\G(2)$ 
\[
\lambda(\tau)\,=\,\left(\frac{\eta(\tau/2)}{\eta(2\tau)}\right)^8.
\]
The following relations will be used below \cite[Chapter 6]{rankin}:
\begin{equation}\label{e4-j}
    E_4\,=\,\frac{J'^2}{(2\pi i)^2J(J-1)}\,,
\end{equation}
\begin{equation}\label{e6-j}
    E_6\,=\,\frac{J'^3}{(2\pi i)^3J^2(J-1)}\,,
\end{equation}
\begin{equation}\label{delta-j}
    \Delta=\frac{-1}{(48\pi^2)^3}\frac{J'^6}{J^4(J-1)^3}\,.
\end{equation}
Let us recall that $\Delta$ does not vanish in $\fH$, while $E_6$ (resp. $E_4$) has a simple zero at $i$ and its $\SL-$ orbit (resp. at $\rho=\exp(2\pi i/3)$). Meanwhile, $J-1$ has a double zero at $i$ and $J$ has a zero of order 3 at $\rho$.

The following propositions will be very useful in the next sections.
\begin{prop}\label{eta-j-diff} We have
\[
    12\,\frac{\eta'(i)}{\eta(i)}\,=\,\frac{3}{7}\frac{E_6''(i)}{E_6'(i)}\,=\,\frac{J'''(i)}{J''(i)}\,=\,3i.
\]
\end{prop}
\begin{proof}
Taking the logarithmic derivative in \eqref{delta-j} yields
\begin{equation}\label{eta44}
24\,\frac{\eta'}{\eta}\,=\,6\,\frac{J''}{J'}-4\,\frac{J'}{J}-3\,\frac{J'}{J-1}.
\end{equation}
Using the expansion of $J$ near $i$
\[
J(\tau)=1+\frac{1}{2}\,J''(i)(\tau-i)^2+\frac{1}{6}\,J'''(i)(\tau-i)^3+\mbox{O}(\tau-i)^4
\]
we get 
\[
\frac{J''(\tau)}{J'(\tau)}=\frac{1}{\tau-i}+\frac{1}{2}\,\frac{J'''(i)}{J''(i)}+\mbox{O}(\tau-i),
\]
\[
\frac{J'(\tau)}{J(\tau)-1}=\frac{2}{\tau-i}+\frac{1}{3}\,\frac{J'''(i)}{J''(i)}+\mbox{O}(\tau-i)
\]
and
\[
\frac{J'(\tau)}{J(\tau)}=\mbox{O}(\tau-i).
\]
Therefore,
\[
24\,\frac{\eta'(\tau)}{\eta(\tau)}=2\,\frac{J'''(i)}{J''(i)}+\mbox{O}(\tau-i),
\]
that is
\[
12\,\frac{\eta'(i)}{\eta(i)}=\frac{J'''(i)}{J''(i)}.
\]
In the meantime, differentiating $J(-1/\tau)=J(\tau)$ trice yields
\[
\frac{6}{\tau^4}\,J'(-1/\tau)-\frac{6}{\tau^5}\,J''(-1/\tau)+\frac{1}{\tau^6}\,J'''(\-1/\tau)\,=\,J'''(\tau).
\]
Since $J'(i)=0$, we get
\[
\frac{J'''(i)}{J''(i)}=3i.
\]
Finally, using 
\[
E_6(\tau)=(\tau-i)E_6'(i)+\frac12 (\tau-i)^2E_6''(i)+\mbox{O}(\tau-i)^3
\]
and
\[
E_6'(\tau)=E_6'(i)+(\tau-i)E_6''(i)+\mbox{O}(\tau-i)^2,
\]
we obtain
\begin{equation} \label{ee6}
    \frac{E_6'(\tau)}{E_6(\tau)}=\frac{1}{\tau-i}+\frac12 \frac{E_6''(i)}{E_6'(i)}+\mbox{O}(\tau-i)^2.
\end{equation}
On the other hand, taking the logarithmic derivative in \eqref{e6-j} yields
\begin{align*}
\frac{E_6'(\tau)}{E_6(\tau)}&=3\frac{J''(\tau)}{J'(\tau)}-2\frac{J'(\tau)}{J(\tau)}-\frac{J'(\tau)}{J(\tau)-1}\\
&=\frac{1}{\tau-i}+\frac76 \frac{J'''(i)}{J''(i)}+\mbox{O}(\tau-i)
\end{align*}
using the expansions of $J''/J'$, $J'/J$ etc. cited in the beginning of this proof. Now, comparing with \eqref{ee6}, we get
\[
\frac{E_6''(i)}{E_6'(i)}=\frac73 \frac{J'''(i)}{J''(i)},
\]
which concludes the proof.
\end{proof}
\begin{prop}\label{eta-e4-diff}
We have
\[
24\frac{\eta'(\rho)}{\eta(\rho)}\,=\,\frac{J^{(4)}(\rho)}{J'''(\rho)}
\,=\,\frac{6}{5}\frac{E_4''(\rho)}{E_4'(\rho)}\,=\,12\,\frac{1+\rho}{1-\rho}.
\]
\end{prop}
\begin{proof}
Write
\[
E_4(\tau)=(\tau-\rho)E_4'(\rho)+\frac{1}{2}(\tau-\rho)^2E_4''(\rho)+\mbox{O}(\tau-\rho)^3,
\]
\[
E_4'(\tau)=E_4'(\rho)+(\tau-\rho)E_4''(\rho)+\mbox{O}(\tau-\rho)^2,
\]
so that
\begin{equation}\label{e44}
\frac{E_4'(\tau)}{E_4(\tau)}\,=\,\frac{1}{\tau-\rho}\,+\,\frac{1}{2}\frac{E_4''(\rho)}{E_4'(\rho)}\,+\,\mbox{O}(\tau-\rho).
\end{equation}
Now write
\[
J(\tau)=\frac{1}{6}J'''(\rho)(\tau-\rho)^3+\frac{1}{24}J^{(4)}(\rho)(\tau-\rho)^4+\mbox{O}(\tau-\rho)^5,
\]
\[
J'(\tau)=\frac{1}{2}J'''(\rho)(\tau-\rho)^2+\frac{1}{6}J^{(4)}(\rho)(\tau-\rho)^3+\mbox{O}(\tau-\rho)^4,
\]
\[
J''(\tau)=J'''(\rho)(\tau-\rho)+\frac{1}{2}J^{(4)}(\rho)(\tau-\rho)^2+\mbox{O}(\tau-\rho)^3.
\]
It follows that
\[
\frac{J''(\tau)}{J'(\tau)}=\frac{2}{\tau-\rho}+\frac{1}{3}\frac{J^{(4)}(\rho)}{J'''(\rho)}+\mbox{O}(\tau-\rho),
\]
\[
\frac{J'(\tau)}{J(\tau)}=\frac{3}{\tau-\rho}+\frac{1}{4}\frac{J^{(4)}(\rho)}{J'''(\rho)}+\mbox{O}(\tau-\rho)
\]
and
\[
\frac{J'(\tau)}{J(\tau)-1}=\mbox{O}(\tau-\rho)^2.
\]
Now, using \eqref{e4-j}, we have
\begin{equation}\label{e55}
\frac{E_4'(\tau)}{E_4(\tau)}=2\frac{J''(\tau)}{J'(\tau)}-\frac{J'(\tau)}{J(\tau)}-\frac{J'(\tau)}{J(\tau)-1}.
\end{equation}
Therefore,
\[
\frac{E_4'(\tau)}{E_4(\tau)}=\frac{1}{\tau-\rho}+\frac{5}{12}\frac{J^{(4)}(\rho)}{J'''(\rho)}+\mbox{O}(\tau-\rho).
\]
Comparing with \eqref{e44}, we obtain
\[
\frac{E_4''(\rho)}{E_4'(\rho)}\,=\,\frac{5}{6}\frac{J^{(4)}(\rho)}{J'''(\rho)}.
\]
On the other hand, using \eqref{eta44}, we get
\[
24\frac{\eta'(\tau)}{\eta(\tau)}\,=\,\frac{J^{(4)}(\rho)}{J'''(\rho)}+\mbox{O}(\tau-\rho),
\]
which proves that
\[
24\frac{\eta'(\rho)}{\eta(\rho}\,=\,\frac{J^{(4)}(\rho)}{J'''(\rho)}
\,=\,\frac{6}{5}\frac{E_4''(\rho)}{E_4'(\rho)}.
\]
Furthermore, differentiating twice the identity 
\[ E_4\left(\frac{-1}{\tau+1}\right)\,=\,(\tau+1)^4E_4(\tau)\]
and taking $\tau=\rho$ yields the last equality in the proposition.
\end{proof}

\section{Level 6 modular differential equations and algebraic systems}
Suppose we are given a $\rho-$equivariant function for a finite index subgroup $\G$ of $\SL$. If $\rho$ is a reducible representation of $\G$, then it can be conjugated to an upper triangular representation, i.e. there exists $\sigma\in\GC$ such that $\rho_1=\sigma\rho\sigma^{-1}$ is upper triangular. Moreover $h_1=\sigma\cdot h$ is $\rho_1-$equivariant and shares the same Schwarz derivative with $h$. Thus, if we are looking for a solution to \eqref{equ2} corresponding to a reducible representation, we may suppose, without loss of generality, that $\rho$ is upper-triangular.
According to Theorem 4.3 in\cite{ramanujan}, a meromorphic function $h$ is $\rho-$equivariant for a triangular representation $\rho$ of $\G$ if and only if the derivative $h'$ is a meromorphic weight 2 modular form for $\G$ with a character. Now, for the Schwarz derivative $\{h,\tau\}$ to be holomorphic,   $h'$ must be nonvanishing where $h$ is holomorphic and, elsewhere, $h$ should have only simple poles, which is equivalent to say that
$h'$ has only double poles with zero residues. Therefore,  if we seek a solution $h$ to \eqref{equ2}, then we have to integrate nonvanishing weight 2 modular forms for $\SL$ with a character and having double poles (if any) with vanishing residues. The characters in question are trivial on the commutator group of $\SL$ which is a level 6 and index 12 in $\SL$, and therefore these modular forms have a $q-$expansion where $q=\exp(2\pi i\tau /6)$.

According to \cite{ramanujan}, a holomorphic weight 2 modular form for $\SL$ with a character must be equal to $c\eta^4$ where $c$ is a constant. If we set
$$ h(\tau)=\int_i^{\tau}\eta^4(\tau)\,d\tau,$$ 
then we have
\[
\{h,\tau\}\,=\,\frac{2\pi^2}{36}\,E_4(\tau).
\]
In other words, $h$ is a solution to \eqref{equ2} with $\disp s=\frac{\pi^2}{6^2}$. It follows that that $y=1/\sqrt{h'}=\eta^{-2}$ is a solution to $\disp y''+\frac{\pi^2}{36}\,y=0$; a differential equation that was first mentioned by Klein in \cite{klein}, and it was Hurwitz who first gave $\eta^{-2}$ as a solution to this equation \cite{hurwitz}.

In order to find other solutions, we should look for weight 2 modular forms with double poles and zero residues. To this end, for each triplet of parameters $(a,b,c)$, we introduce the following algebraic system $E_{a,b,c}^n$ of $n$ equations in $n$ variables $x_i$:
\begin{equation}\label{system}
\frac{a}{x_i}\,+\,\frac{b}{x_i-1}\,+\,\sum_{j\neq i}\,\frac{c}{x_j-x_i}\,=\, 0\,,\ \ 1\leq i\leq n.
\end{equation}
Notice that for $\alpha\neq 0$, the system $E^n_{a,b,c}$ is equivalent to the system $E^n_{\alpha a,\alpha b,\alpha c}$.
According to \cite[Theorem 6.2]{ramanujan}, if $a$, $b$ and $c$ are positive real numbers, then the system $E_{a,b,c}^n$ has a solution  in $(0,1)^n$.
Let $(x_i)_{1\leq i\leq n}$ be a solution to the algebraic system $E^n_{4,3,12}$ and set
\[
f_n\,=\,\frac{\eta^4}{\prod_{i=1}^n\,(J-x_i)^2}.
\]
Also, write $x_i=J(w_i)$, $w_i\in\fH$. Then, as $0<x_i<1$, the $w_i$'s are not elliptic fixed points and $f_n$ is a weight 2 modular form with a character and has a double pole at each $w_i$, $1\leq i\leq n$, and is holomorphic elsewhere. Moreover, the fact that the $x_i$'s satisfy the system $E_{4,3,12}^n$ is equivalent to the vanishing of the residues of $f_n$ at each $w_i$. One of the main results in \cite{ramanujan} is that $\disp h_n(\tau)=\int_i^{\tau}\,f_n(z)\,dz$ is a solution to \eqref{equ2}
with $\disp s=\pi^2\left(\frac{12n+1}{6}\right)^2$, while the solutions to \eqref{equ1} are given by $\disp y_1=\eta^{-2}\prod_{i=1}^n\,(J-x_i)$ and $y_2=h_ny_1$; generalizing Hurwitz's solution when $n=0$. This solves the modular differential equations \eqref{equ1} and \eqref{equ2} with 
$s=\pi^2(m/6)^2$ with $m\equiv 1 \mod 12$. We now focus on finding the solutions for the remaining residues classes modulo 12, that is, when $m\equiv 5$, 7 or 11 $\mod 12$. The idea is to allow  the nonvanishing weight
2 modular forms to have double poles at elliptic points.
\begin{thm}
Let $n\in{\mathbb N}$ and $(x_i)_{1\leq i\leq n}\in(0,1)^n$ be a solution to the algebraic system $E_{4,9,12}^n$ and let $w_i\in\fH$ such that $x_i=J(w_i) $. Then
\[
f_n\,=\, \frac{\eta^4}{(J-1)\prod_{i=1}^n(J-x_i)^2}\,=\, \frac{\eta^{28}}{E_6^2\prod_{i=1}^n(J-x_i)^2}
\]
is a nonvanishing weight 2 modular form with double poles at $i$ and at each $w_i$ with zero residues. Moreover, if $\disp h_n(\tau)=\int_i^{\tau}f_n(z)dz$, then
\[
\{h_n,\tau\}\,=\,2\pi^2\frac{(12n+7)^2}{6^2}\,E_4(\tau).
\]
\end{thm}
\begin{proof}
As $i$ is a double zero of $J-1$ and each $w_i$ is not in the $\SL-$orbit of $i$, it is clear that $i$ is a double pole of $f_n$.
Write $\disp g_n=\eta^4/\prod_{i=1}^n(J-x_i)^2$ so that $f_n=g_n/(J-1)$. Also write
\[
J(\tau)-1=\frac{1}{2}J''(i)(\tau-i)^2+\frac{1}{6}J'''(i)(\tau-i)^3+\mbox{O}(\tau-i)^4.
\]
Then
\[
\frac{J''(i)}{2}\frac{(\tau-i)^2}{J(\tau)-1}=1-\frac{1}{3}\frac{J'''(i)}{J''(i)}(\tau-i)+\mbox{O}(\tau-i)^2.
\]
It follows that
\[
\frac{d}{d\tau}((\tau-i)^2f_n(\tau))=\frac{d}{d\tau}\frac{g_n(\tau)(\tau-i)^2}{J(\tau)-1}
=g_n'(\tau)\frac{(\tau-1)^2}{J(\tau)-1}+g_n(\tau)\frac{d}{d\tau}\frac{(\tau-i)^2}{J(\tau)-1}.
\]
Therefore,
\[
\mbox{Res}(f_n,i)=\lim_{\tau\rightarrow i}\frac{d}{d\tau}((\tau-i)^2f_n(\tau))
=
\frac{2g_n(i)}{J''(i)}\left(\frac{g_n'(i)}{g_n(i)}-\frac{J'''(i)}{3J''(i)}\right).
\]
In the meantime, taking the logarithmic derivative of $g_n$ yields
\[
\frac{g_n'(i)}{g_n(i)}=4\frac{\eta'(i)}{\eta(i)}-\sum_{i=1}^n\frac{2J'(i)}{J(i)-x_i}=4\frac{\eta'(i)}{\eta(i)}.
\]
Hence, using \propref{eta-j-diff},
\[
\mbox{Res}(f_n,i)=\frac{2g_n(i)}{J''(i)}\left(4\frac{\eta'(i)}{\eta(i)}-\frac{J'''(i)}{3J''(i)}\right)=0.
\]

Now, fix $i$, $1\leq i\leq n$, and write $f_n(\tau)=\phi_n(\tau)/(J(\tau)-J(w_i))^2$.
A similar calculation as above shows that
\[
\mbox{Res}(f_n,w_i)=\frac{\phi_n'(w_i)}{J'(w_i)^2}-\phi_n(w_i)\frac{J''(w_i)}{J'(w_i)^3}=\frac{\phi_n(w_i)}{J'(w_i)^2}\left(\frac{\phi_n'(w_i)}{\phi_n(w_i)}-\frac{J''(w_i)}{J'(w_i)}  \right).
\]
Meanwhile,
\[
\frac{\phi_n'}{\phi_n}=\frac{4\eta'}{\eta}-\frac{J'}{J-1}-\sum_{j\neq i}\frac{2J'}{J-J(w_j)}.
\]
Therefore
\[
\mbox{Res}(f_n,w_i)=\frac{\phi_n(w_i)}{J'(w_i)^2}\left(\frac{4\eta'(w_i)}{\eta(w_i)}-\frac{J'(w_i)}{J(w_i)-1}-\sum_{j\neq i}\frac{2J'(w_i)}{J(w_i)-J(w_j)}-\frac{J''(w_i)}{J'(w_i)}  \right).
\]
Using \eqref{eta44}, we have
\[
\frac{4\eta'}{\eta}-\frac{J''}{J'}=\frac{J'}{6}\left(\frac{3}{1-J}-\frac{4}{J}\right),
\]
and hence 
\[
\mbox{Res}(f_n,w_i)=\frac{-\phi_n(w_i)}{6J'(w_i)}\left(\frac{4}{J(w_i)}+\frac{9}{J(w_i)-1}+\sum_{j\neq i}\frac{12}{J(w_i)-J(w_j)}
\right)=0
\]
because the $x_i=J(w_i)$ were chosen to be a solution to the algebraic system $E^n_{4,9,12}$. Therefore $f_n$ has only double poles with zero residues and is nonvanishing elsewhere. Thus its integral $h_n$ has only simple poles and it is locally univalent elsewhere. It follows that the Schwarz derivative $\{h_n,\tau\}$ is a holomorphic weight 4 modular form for $\SL$ and hence a scalar multiple of $E_4$. Finally, notice that the leading coefficient of the $q-$expansion of $f_n$ is $q^{\frac{1}{6}+1+2n}=q^{\frac{12n+7}{6}}$ and consequently the leading coefficient of
$\{h_n,\tau\}=(f_n'/f_n)'-\frac{1}{2}(f_n'/f_n)^2$ is easily seen to be 
$\disp 2\pi^2\frac{(12n+7)^2}{6^2}$.

\end{proof}
We now seek a similar solution but with a double at the other elliptic fixed points, namely $\tau=\rho$.
\begin{thm}
Let $n\in{\mathbb N}$ and $(x_i)_{1\leq i\leq n}\in(0,1)^n$ be a solution to the algebraic system $E_{8,3,12}^n$ and let $w_i\in\fH$ such that $x_i=J(w_i) $. Then
\[
f_n\,=\, \frac{\eta^{20}}{E_4^2\prod_{i=1}^n(J-x_i)^2}
\]
is a nonvanishing weight 2 modular forms with double poles at $\rho$ and at each $w_i$ with zero residues. Moreover, if $\disp h_n(\tau)=\int_i^{\tau}f_n(z)dz$, then
\[
\{h_n,\tau\}\,=\,2\pi^2\frac{(12n+5)^2}{6^2}\,E_4(\tau).
\]
\end{thm}
\begin{proof}
Write $f_n=\psi_n/E_4^2$ and
\[
 E_4(\tau)=(\tau-\rho)E_4'(\rho)+ \frac{1}{2}(\tau-\rho)^2E_4''(\rho)+\mbox{O}(\tau-\rho)^3
\]
so that
\[
E_4^2(\tau)=(\tau-\rho)^2E_4'^2(\rho)\left(1+(\tau-\rho)\frac{E_4''(\rho)}{E_4'(\rho)}+\mbox{O}(\tau-\rho)^2   \right).
\]
Hence,
\[
\frac{(\tau-\rho)^2\psi_n(\tau)}{E_4^2(\tau)}=\frac{\psi_n(\tau)}{E_4'^2(\rho)}
\left(1-(\tau-\rho)\frac{E_4''(\rho)}{E_4'(\rho)}+\mbox{O}(\tau-\rho)^2   \right).
\]
It follows that
\begin{align*}
\lim_{\tau\rightarrow \rho}\,
\frac{d}{d\tau}\,\frac{(\tau-\rho)^2\psi_n(\tau)}{E_4^2(\tau)}
&=\frac{\psi_n'(\rho)}{E_4'^2(\rho)}-\frac{\psi_n(\rho)E_4''(\rho)}{E_4'^3(\rho)}\\
&=\frac{\psi_n(\rho)}{E_4'^2(\rho)}\left(\frac{\psi_n'(\rho)}{\psi_n(\rho)}-\frac{E_4''(\rho)}{E_4'(\rho)} \right)\\
&= \frac{\psi_n(\rho)}{E_4'^2(\rho)}\left(20\frac{\eta'(\rho)}{\eta(\rho)}-\sum_{i=1}^n\frac{2J'(\rho)}{J(\rho)-x_i)}-
\frac{E_4''(\rho)}{E_4'(\rho)} \right)\\
&=\frac{\psi_n(\rho)}{E_4'^2(\rho)}\left(20\frac{\eta'(\rho)}{\eta(\rho)}-
\frac{E_4''(\rho)}{E_4'(\rho)} \right)\\
&=0.  
\end{align*}
The last equality follows from \propref{eta-e4-diff}.
Therefore, the residue of $f_n$ at the double pole $\rho$ is zero. In a similar manner to the previous theorem and using both \eqref{eta44} and \eqref{e55}, it is easily seen that the residues of $f_n$ at each $w_i$ is precisely zero  because the $x_i$'s satisfy the algebraic system $E^n_{8,3,12}$.
Finally, the leading coefficient of the $q-$expansion of $f_n$ is $\disp q^{\frac{12n+5}{6}}$ and thus the leading coefficient of $\{h_n,\tau\}$ is
$\disp 2\pi^2\frac{(12n+5)^2}{6^2}$.

\end{proof}

Finally, we seek a solution which has both elliptic points $i$ and $\rho$ as double poles.
\begin{thm}
Let $n\in{\mathbb N}$ and $(x_i)_{1\leq i\leq n}\in(0,1)^n$ be a solution to the algebraic system $E_{8,9,12}^n$ and let $w_i\in\fH$ such that $x_i=J(w_i) $. Then
\[
f_n\,=\, \frac{\eta^{20}}{E_4^2(J-1)\prod_{i=1}^n(J-x_i)^2}\,=\,
\frac{\eta^{44}}{E_4^2E_6^2\prod_{i=1}^n(J-x_i)^2}
\]
is a nonvanishing weight 2 modular forms with double poles at $i$,  $\rho$ and at each $w_i$ with zero residues. Moreover, if $\disp h_n(\tau)=\int_i^{\tau}f_n(z)dz$, then
\[
\{h_n,\tau\}\,=\,2\pi^2\frac{(12n+11)^2}{6^2}\,E_4(\tau).
\]
\end{thm}
\begin{proof}
This can be shown in the same way as the previous two theorems  with the use of both \propref{eta-j-diff} and \propref{eta-e4-diff}. At the same time, the exponent of $q$ in the leading coefficient of $f_n$ is $\disp \frac56 +1+12n=\frac{12n+11}{6}$.
\end{proof}

\section{Modular solutions and algebraic systems}
We have mentioned  that according to \cite{forum}, the Schwarzian equation	\eqref{equ2} has solutions that are modular functions  if and only if $s=\pi^2n^2/m^2$ with $m$ and $n$ being positive integers such that $2\leq m\leq 5$ and $\gcd(m,n)=1$. For each such pair $(m,n)$, the invariance group for the modular solution $h$ is $\G(m)$ and $n$ is the ramification index above $\infty$ in the covering
$h:X(m)\longrightarrow \PC$. Here $X(m)=X(\G(m))$. A key fact about the groups $\G(m)$ for $2\leq m\leq 5$ is that they are the only principal congruence groups that are genus 0 and torsion-free. In this section, we will establish that these modular solutions are also attached to an algebraic system in the same way the solutions in the previous section were.

Let $m\in\{2,3,4,5\}$ and let $t$ be a Hauptmodul of $\G(m)$. Choose $t$ so that its Fourier expansion has the shape
\[
t(\tau)=\frac{1}{q}+\sum_{i\geq 0}\,a_iq^i\,,\ \ \ q=e^\frac{2\pi i\tau}{m}.
\]
Since the Hauptmodul $t$ takes its values only once and $\G(m)$ has  no elliptic elements, then according to \cite{mathann}, $\{t,\tau\}$ is a holomorphic weight 4 modular form for the normalizer of $\G(m)$ in $\SL$ which is $\SL$ itself , and thus it is a scalar multiple of $E_4$. From the $q-$expansion of $t$, it is clear that
\[
\{t,\tau\}\,=\,\frac{\pi^2}{m^2}\,E_4(\tau).
\]
Now let $n\geq 2$ be a integer coprime with $m$. According to \cite{forum}, there exists a modular function $h$ for $\G(m)$ solution to $\{h,\tau\}=2\pi^2(n/m)^2E_4(\tau)$. As $n$ is the ramification index of $h$ at $\infty$, we can write
\[
h(\tau)=q^n+\mbox{o}(q^n)\,\,\ \ q=e^\frac{2\pi i\tau}{m}.
\]
Now, suppose that the poles of $h$ are given by the set $\{w_1,w_2,\ldots w_a, s_1,\ldots s_b\}$, where, for $1\leq i\leq a$, $w_i\in\fH$ (if any)  and the $s_j$, $1\leq j\leq b$, are among the cusps of $\G(m)$. Then the degree $d$ of the covering
$h:X(m)\longrightarrow {\mathbb P}^1(\BC)$ satisfies
\begin{equation}\label{danb}
d\,=\,a+nb.
\end{equation}
We also consider the modular function $f=t'/h'$ for $\G(m)$. Since $h'$ can have only double poles at the $w_i$'s and it is nonvanishing elsewhere in $\fH$, we see that $f$ is holomorphic in $\fH$. Therefore, for some polynomials $P$ and $Q$, we have $f=P(t)/Q(t)$. Moreover, as the Hauptmodul $t$ has a pole at $\infty$, it is holomorphic on $\fH$ and $t'$ does not vanish on $\fH$ (because $t$ is a Hauptmodul and $\G(m)$ has no elliptic element, or because the Schwarz derivative of $t$ is holomorphic as we have seen above). It follows that each $w_i$, $1\leq i\leq a$, is a zero of order 2 of $f$. In the meantime, the behaviour of $f$ at the cusps is as follows:
\begin{itemize}
    \item Near each $s_i$, $1\leq i\leq b$: we have, for some constants $\alpha$, $\beta$ and $\gamma$, $h'(\tau)=\alpha/q^n + \ldots$ and $t'(\tau)=\beta q+\ldots$ because $t$ has a pole at $\infty$ and thus it is holomorphic at any other cusp. Therefore, $$f(\tau)=\gamma q^{n+1}+\ldots.$$
    \item Near $\infty$: \[\displaystyle f(\tau)=\frac{\alpha/q+\ldots}{\beta q^n+\ldots}=\frac{\gamma}{q^{n+1}}+\ldots.\]
    \item Near each cusp $s\notin\{s_1,\ldots,s_b,\infty\}$: we have \\ \[f(\tau)=\frac{\alpha q+\ldots}{\beta q^n+\ldots}=\frac{\gamma}{q^{n-1}}+\ldots.\]
\end{itemize}
Therefore, we have
\[
P(t)=\prod_{i=1}^a\,(t-t(w_i))^2\,\prod_{i=1}^b\,(t-t(s_i))^{n+1}
\]
and
\[
Q(t)=\prod_{s\notin\{s_1,\ldots,s_b,\infty\}}(t-t(s))^{n-1}.
\]
Furthermore, comparing the order of $\infty$ in $h'/t'=Q(t)/P(t)$ yields
\[
(n+1)=2a+(n+1)b-(n-1)(\nu_{\infty}-(b+1)),
\]
where $\nu_{\infty}$ is the number of inequivalent cusps for $\G(m)$. Hence, using \eqref{danb}, we get
\[
2d-2=(n-1)\nu_{\infty}.
\]
 Notice that this is simply the Riemann-Hurwitz formula for the  covering $h:X(m)\longrightarrow {\mathbb P}^1(\BC)$.

We can have a more precise information on $a$ and $b$ for a given level $m$.
 \begin{prop}
     With the notation as above, for each positive integer $n$, we have
     \begin{enumerate}
         \item If $m=2$, then 
         \[
         (a,b)\in\left\{\left(\frac{3n-1}{2},0\right), \left(\frac{n-1}{2},1\right)\right\}.
         \]
         \item If $m=3$, then
         \[
         (a,b)\in\{(2n-1,0),(n-1,1)\}.
         \]
         \item If $m=4$, then
         \[
         (a,b)\in\{(3n-2,0),(2n-2,1),(n-2,2)\}.
         \]
         \item Finally, if $m=5$, then
         \[
         (a,b)\in\{(n(6-k)-5,k),\, 0\leq k\leq 5\}.
         \]
     \end{enumerate}
 \end{prop}
	\begin{proof}
	    If $m=2$, there are 3 inequivalent cusps and the Riemann-Hurwitz formula reads $2d=3n-1$ which implies that $n$ is odd. Since $a+nb=d$, we have then $2a+1=n(3-2b)$. It follows that either $b=0$ which gives $a=(3n-1)/2$ or $b=1$ which corresponds to $a=(n-1)/2$.
Similarly, for $m=3$, $\nu_{\infty}=4$ and then
     $a+nb=2n-1$ which we can rewrite as $a+1=n(2-b)$. It follows that $b\in\{0,1\}$ and $(a,b)\in\{(2n-1,0),(n-1,1)\}$. The other cases follow similarly knowing that for $n\geq 3$, $\disp \nu_{\infty}=\frac12 n^2\prod_{p|n, p\, prime}(1-1/p^2)$.
	\end{proof}
 Finally, the solution to $\{h,\tau\}=2\pi^2(n/m)^2 E_4$ is thus obtained by integrating the weight 2 modular form
 $h'=t'Q(t)/P(t)$  by choosing the adequate pair $(a,b)$ given in the above proposition. Clearly $t'Q(t)$ does not vanish on $\fH$ and the poles in $\fH$ are the $w_i$'s which should have a zero residue. Hence, $x_i=t(w_i)$, $1\leq i\leq a$, are a solution to a system of type \eqref{system}. 
 
 Let us illustrate this construction in the case $m=2$. The group $\G(2)$ has 3 cusps, namely 0, 1 and $\infty$. We take $t=1/\la$ which sends the triple $(0,1,\infty)$ to the triple $(1,0,\infty)$.

 \underline{\bf Case 1:} If $a=(3n-1)/2$ and $b=0$, then 
 \begin{equation}\label{case1}
 h'=\frac{t't^{n-1}(t-1)^{n-1}}{\prod_{j=1}^a (t-t(w_j))^2}.
 \end{equation}
 Fix $i\in\{1\ldots a\}$, and write $h'=g/(t-t(w_i))^2$. Then 
 \[
 \mbox{Res}(h',w_i)=\frac{g(w_i)}{t'(w_i)^2}\left(\frac{g'(w_i)}{g(w_i)}-\frac{t''(w_i)}{t'(w_i)}\right).
 \]
Meanwhile,
\[
\frac{g'}{g}=\frac{t''}{t'}+(n-1)\frac{t'}{t}+(n-1)\frac{t'}{t-1}-\sum_{j\neq i}\,\frac{2t'}{t-t(w_j)}.
\]
 It follows that 
 \[
 \mbox{Res}(h',w_i)=\frac{g(w_i)}{t'(w_i)}\left(\frac{n-1}{t(w_i)}-\frac{n-1}{t(w_i)-1}-\sum_{j\neq i}\,\frac{2}{t(w_i)-t(w_j)}\right).
 \]
 Thus, if we set $x_i=t(w_i)$, then $(x_i)_{1\leq i\leq a}$ is a solution to the algebraic system
$E^a_{n-1,1-n,-2}$.

\underline{\bf Case 2:} If $a=(n-1)/2$ and $b=1$, then we have two sub-cases depending on whether we take the cusp 0 or the cusp 1 in the polynomial $P(t)$. We have the two possibilities for $h'$:
\begin{equation}\label{case2}
h'_1=\frac{t't^{n-1}}{(t-1)^{n+1}\prod_{j=1}^a (t-x_j)^2}\ \ \ \mbox{and}\ \ \ 
h'_2=\frac{t'(t-1)^{n-1}}{t^{n+1}\prod_{j=1}^a (t-x_j)^2},
\end{equation}
where the $(x_i)_{1\leq i\leq a}$ are solution to $E^a_{1-n,n+1,2}$ and $E^a_{n+1,1-n,2}$ respectively.
Notice that both functions have $q^n$ as a leading term and their integrals are solution to $\{h,\tau\}=\pi^2(n/2)^2E_4(\tau)$. This mean that $h_1$ and $h_2$ are linear fraction of one another.

\begin{example}{\rm  For $m=2$ and $n=3$, we have three solutions for $h'$:
\begin{enumerate}
    \item With one pole in $\fH$ and one pole at the cusp 0:
\[
h'_1=\frac{t't^2}{(t-x)^2(t-1)^4}.
\]
The residue at the pole in $\fH$ is zero lead to $x=-1$. Thus, we have the solution
\[
h_1=\frac{-1}{6}\frac{2t-1}{(t-1)^3(t+1)}.
\]
\item With one pole in $\fH$ and one pole at the cusp 1:
\[
h'_2=\frac{t'(t-1)^2}{(t-x)^2t^4}.
\]
The residue at the pole in $\fH$ vanishes when $x=2$. The primitive is given by
\[
h_2=\frac{-1}{6}\frac{2t-1}{t^3(t-2)}.
\]
\item With four poles in $\fH$ and none at the cusps:
\[
h'_3=\frac{t' t^2(t-1)^2}{\prod_{i=1}^4(t-x_i)},
\]
where $x_1,\ldots,x_4$ are solutions to
\[
\frac{1}{x_i}+\frac{1}{x_i-1}-\sum_{j\neq i}\,\frac{1}{x_i-x_j}\,,\ \ 1\leq i\leq 4.
\]
This algebraic system has as solutions (up to a permutation):
\[
\frac{1-\sqrt{3}}{2}\pm \left(\frac34\right)^{\frac14}\ ,\ \ \ \frac{1+\sqrt{3}}{2}\pm i\left(\frac34\right)^{\frac14}.
\]
These solutions are the roots of the irreducible polynomial 
\[
P(x)=x^4-2x^3+4x-2.
\]
Therefore we can write
\[
h'_3=\frac{t't^2(t-1)^2}{(t^4-2t^3+4t-2)^2}\ ,\ \ \ h_3=\frac{1}{12}\frac{t^3(t-2)}{t^4-2t^3+4t-2}.
\]

\end{enumerate}
 }
\end{example}
\begin{remark}{\rm
We expect that for each $ 1\leq m\leq 5 $, every choice of the pair $(a,b)$ gives arise to a solution of the  same Schwarz differential equation $\{h,\tau\}=2\pi^2(n/m)^2 E_4$, and hence these solutions should be linear fractions of each others.
This is illustrated in the case of  $m=2$ and $n=3$, where it can be easily checked that
\[
h_2\,=\,\frac{h_1}{6h_1+1}\ \mbox{ and }\  h_3\,=\,\frac{6h_1+1}{-72h_1+12}.
\]
}
\end{remark}
\begin{remark}{\rm
The lambda function and its derivative can be expressed in terms of Jacobi theta functions. Indeed, according to \cite[Chapter 7]{rankin}, we have
\[
\lambda=\frac{\theta_2^4}{\theta_3^4}\ ,\ \ \lambda'=i\pi \theta_4^4\lambda=i\pi \frac{\theta_2^4\theta_4^4}{\theta_3^4}.
\]
Hence
\[
t=\frac{\theta_3^4}{\theta_2^4}\ ,\ \ t'=-i\pi \frac{\theta_3^4\theta_4^4}{\theta_2^4}.
\]
Therefore, we can see that in the general case for the level 2, the derivatives in \eqref{case1} and in \eqref{case2} are readily squares since $n$ is odd. This allows to easily write down a square root of $h'$ whose reciprocal is a solution to the modular differential equation $y''+\pi^2(n/2)^2E_4\,y=0$.
}
\end{remark}



\begin{thebibliography}{aaaa}
 \bibitem{yang} Z. Chen; C-S. Lin; Y. Yang. Modular ordinary differential equations on $\SL$ of third order and applications. SIGMA Symmetry Integrability Geom. Methods Appl. 18 (2022), Paper No. 013.
		\bibitem{structure} A. Elbasraoui; A. Sebbar. Equivariant forms: Structure and geometry. Canad. Math. Bull. Vol. {\bf 56} (3), (2013) 520--533.
		\bibitem{rational} A. Elbasraoui; A. Sebbar. Rational equivariant forms. Int. J. Number Th. 08  No. 4(2012), 963--981.
	\bibitem{fr-ma} C. Franc; G. Mason. Hypergeometric series, modular linear differential equations and vector-valued modular forms. The Ramanujan J. 41, 233--267 (2016). 	
  \bibitem{grabner} P. J. Grabner. Quasimodular forms as solutions of modular differential equations.  Int. J. Number Theory 16 (2020), no. 10, 2233--2274.
		\bibitem{hurwitz} A. Hurwitz, Adolf: Ueber die Differentialgleichungen dritter Ordnung, welchen
		\bibitem{ka-ko} M. Kaneko; M. Koike, On modular forms arising from a differential equation of hypergeometric type. Ramanujan J. 7(2003), no. 1--3, 145--164.
		\bibitem{ka-et-al} M. Kaneko; K. Nagatomo; Y. Sakai. The third order modular linear differential equations. J. Algebra 485 (2017), 332--352.
		\bibitem{ka-za} M. Kaneko; D. Zagier. Supersingular j-invariants, hypergeometric series, and Atkin's orthogonal polynomials. Computational perspectives on number theory (Chicago, IL, 1995), 97--126, AMS/IP Stud. Adv. Math., 7, Amer. Math. Soc., Providence, RI, 1998.
		\bibitem{klein} F. Klein,  Ueber Multiplicatorgleichungen. (German) Math. Ann. 15 (1879), no. 1, 86--88.
		\bibitem{milas} A. Milas, Ramanujan’s “Lost Notebook” and the Virasoro algebra. Comm. Math. Phys. 251(2004), no. 3, 657--678.
		\bibitem{mukhi} S. Mathur, S. Mukhi, and A. Sen, On the classification of rational conformal field theories. Phys. Lett. B 213(1988), no. 3, 303--308.
		\bibitem{mathann} J. McKay; A. Sebbar.  Fuchsian groups, automorphic functions
		and Schwarzians. Math. Ann. 318 (2), (2000) 255--275.
		Ramanujan J (2008) 17: 405--427.
		\bibitem{nakaya} T. Nakaya. On modular solutions of certain modular differential equation and supersingular polynomials. Ramanujan J. 48 (2019), no. 1, 13–-20. 
		\bibitem{nehari} Z. Nehari, (1949), The Schwarzian derivative and schlicht functions, Bulletin of the American Mathematical Society, 55 (1949) 545--551.
		\bibitem{rankin} R. Rankin. Modular Forms and Functions, Cambridge Univ. Press, Cambridge, 1977.
		\bibitem{s-g} G. Sansone, J. Gerretsen. Lectures on the theory of functions of a complex variable. II,	Geometric theory. Wolters--Noordhoff Publishing, Groningen 1969.
		\bibitem{forum} H. Saber; A. Sebbar. Automorphic Schwarzian equations.  Forum Math. 32 (2020), no. 6, 1621--1636. 
		\bibitem{ramanujan} H. Saber; A. Sebbar. Automorphic Schwarzian equations and integrals of weight 2 forms. The Ramanujan J.  57 (2022), no. 2, 551--568. 
		\bibitem{jmaa} H. Saber; A. Sebbar. Equivariant solutions to modular Schwarzian equations, J. Math. Anal. Appl., 508 (2022), no. 2,  Paper No. 125887. 
		\bibitem{critical} A. Sebbar; H. Saber. On the critical points of modular forms.  J. Number Theory 132 (2012), no. 8, 1780--1787.
		\bibitem{vvmf} A. Sebbar; H. Saber. Equivariant functions and vector-valued modular forms. Int. J. Number Theory 10 (2014), no. 4, 949--954.
		\bibitem{kyushu} A. Sebbar; H. Saber. On the existence of vector-valued automorphic forms. Kyushu J. Math. 71 (2017), no. 2, 271--285.
        \bibitem{baus} H. Saber; A. Sebbar.    On the modularity of solutions to certain differential equations of hypergeometric type. Bull. Aust. Math. Soc. 105 (2022), no. 3, 385--391.
  \bibitem{vdp} B. Van der Pol; On a non-linear partial differential equation satisfied by the logarithm of the Jacobian theta-functions, with arithmetical applications. I, II. Nederl. Akad. Wetensch. Proc. Ser. A. 54 = Indagationes Math. 13, (1951). 261--271, 272--284.
	\end{thebibliography}
	\end{document}